\documentclass{article}

\usepackage{arxiv}

\usepackage[utf8]{inputenc} 
\usepackage[T1]{fontenc}    
\usepackage{hyperref}       
\usepackage{url}            
\usepackage{booktabs}       
\usepackage{amsfonts}       
\usepackage{nicefrac}       
\usepackage{microtype}      
\usepackage{lipsum}		
\usepackage{graphicx}
\usepackage{natbib}
\usepackage{doi}
\usepackage{amsmath}
\usepackage{algorithm}
\usepackage{algpseudocode}
\usepackage{caption}
\usepackage{subcaption}

\title{Differential Dynamic Programming with Stagewise Equality and Inequality Constraints using Interior Point Method}

\author{ 
    {\hspace{1mm}Siddharth Prabhu}\\
	Department of Chemical \\ and Biomolecular Engineering\\
	Lehigh University\\
	Bethlehem, PA 18015 \\
	\texttt{scp220@lehigh.edu} \\
	\And
    {\hspace{1mm}Srinivas Rangarajan} \\
	Department of Chemical \\ and Biomolecular Engineering\\
	Lehigh University\\
	Bethlehem, PA 18015 \\
	\texttt{srr516@lehigh.edu} \\
        \And
    {\hspace{1mm}Mayuresh Kothare} \\
	Department of Chemical \\ and Biomolecular Engineering\\
	Lehigh University\\
	Bethlehem, PA 18015 \\
	\texttt{mvk2@lehigh.edu} \\
}

\renewcommand{\shorttitle}{}

\begin{document}
\maketitle

\begin{abstract}
	Differential Dynamic Programming (DDP) is one of the indirect methods for solving an optimal control problem. Several extensions to DDP have been proposed to add stagewise state and control constraints, which can mainly be classified as augmented lagrangian methods, active set methods, and barrier methods. In this paper, we use an interior point method, which is a type of barrier method, to incorporate arbitrary stagewise equality and inequality state and control constraints. We also provide explicit update formulas for all the involved variables. Finally, we apply this algorithm to example systems such as the inverted pendulum, a continuously stirred tank reactor, car parking, and obstacle avoidance.
\end{abstract}


\section{Introduction}
One of the popular algorithms to solve an optimal control problem in recent times, is called Differential Dynamic Programming (DDP) \cite{MAYNE}. It is an indirect method that finds the optimal control law by minimizing the quadratic approximation of the value function. A closely related algorithm is the Iterative Linear Quadratic Regulator (iLQR) \cite{Li2004IterativeLQ}, which skips the second-order approximation terms of the dynamic system.  

DDP in its original form does not admit state and control constraints. There have been several extensions to overcome this drawback, and most of them, if not all, fall into three major classes of solution techniques - augmented lagrangian methods, active set methods, and barrier methods. All of these methods essentially solve a two-layer optimization problem. In the augmented lagrangian method \cite{ALTRO, jallet2022constrained, aoyama2020constrained, 8206457}, the outer layer transforms the constrained optimization problem into an unconstrained optimization problem by incorporating the constraints into the objective function with an appropriate penalty. This unconstrained problem is then solved using DDP but with an augmented cost function. The outer loop also updates the penalties based on the infeasibility of each constraint. Augmented lagrangian methods are known to cause numerical issues in the hessian as the penalties increase. In active set methods \cite{7989086, 1104123, multireservoir}, the outer layer guesses an active set, thereby converting an inequality constraint optimization problem into an equality constraint optimization problem. In the inner loop, an equality-constrained problem is solved using DDP. The outer layer also updates the active set based on the solution of the inner loop. One of the major challenges of these methods is identifying or modifying the active set, which can lead to combinatorial complexity in the worst case. In barrier methods \cite{8671755, 8317745, IPDDP}, the outer loop adds inequality constraints in the objective function using barrier functions. The inner loop solves a sequence of optimization problems using DDP with a decreasing value of the barrier parameter. 

In this paper, we build upon the interior point differential dynamic programming proposed in \cite{IPDDP}. Specifically, we extend the algorithm by adding equality constraints along with the inequality constraints that were already proposed. We also provide explicit equations for updating the Lagrange multipliers and the slack variables corresponding to the inequality constraints, which were not provided previously. These explicit update rules are obtained by solving an inexact Newtons method proposed in \cite{FREY20206522}

\section{Differential Dynamic Programming}
\subsection{Preliminaries}
Consider a discrete-time dynamical system 
\begin{equation}
\label{eqn:dyn}
    x_{k + 1} = f(x_k, u_k)
\end{equation}

where $x_k \in \mathbb{R}^n$ and $u_k \in \mathbb{R}^m$ are the states and control inputs at time $k$. The function $f : \mathbb{R}^n \times \mathbb{R}^m \rightarrow \mathbb{R}^n$ describes the evolution of the states to time $k + 1$ given the states and control inputs at time $k$. Consider a finite time optimal control problem starting at initial state $x_0$

\begin{align}
\label{eqn:opt-cont}
\begin{split}
    J^*(\textit{X}, \textit{U}) & = \min _{\textit{U}} \sum _{k = 0}^{N - 1} l(x_k, u_k) + l_f(x_N) \\
    \text{subject to} & \\
    & \quad x_{k + 1} = f(x_k, u_k)    
\end{split}
\end{align}

where the scalar-valued functions $l, l_f, J$ denote the running cost, terminal cost, and total cost respectively. $\textit{X} = (x_0, x_1, ..., x_N)$ and $ \textit{U} = (u_0, u_1, ..., u_{N - 1})$ are the sequence of state and control inputs over the control horizon $N$. We can solve this problem using dynamic programming. If we define the optimal value function at time $k$ as 

\begin{equation}
\label{eqn:val-fun}
    V_k(x_k) = \min _{u_k} \left[ l(x_k, u_k) + V_{k + 1}(x_{k + 1}) \right]
\end{equation}

then, starting from $V_N(x_N) = l_f(x_N)$, the solution to the finite time optimal control problem in \ref{eqn:opt-cont} boils down to finding $V_0$. At every time step $k$, differential dynamic programming solves the optimization problem in \ref{eqn:val-fun} using a quadratic approximation of the value function. The value function is approximated around the states obtained by integrating \ref{eqn:dyn} for given control inputs. Let $Q(x_k + \delta x_k, u_K + \delta u_k)$ be the second order Taylor series approximation of \ref{eqn:val-fun} around the point $x_k$ and $u_k$ then equation \ref{eqn:val-fun}, after dropping the subscript $k$ for simplicity, can be written as 

\begin{equation}
    V_k(x) = \min _{\delta u} Q +  \begin{bmatrix}
    Q_x \\ Q_u \end{bmatrix}^T \begin{bmatrix} \delta x \\ \delta u\end{bmatrix} + \frac{1}{2} \begin{bmatrix} \delta x \\ \delta u \end{bmatrix}^T \begin{bmatrix} Q_{xx} & Q_{xu} \\ Q_{ux} & Q_{uu} \end{bmatrix} \begin{bmatrix} \delta x \\ \delta u \end{bmatrix}
\end{equation}
where 
\begin{align}
\label{eqn:qder}
\begin{split}
    & Q = l + V\\
    & Q_x = l_x + f_x^TV_x \\ 
    & Q_u = l_u + f_u^TV_x \\ 
    & Q_{xx} = l_{xx} + f_x^TV_{xx}f_x + V_x f_{xx} \\ 
    & Q_{xu} = l_{xu} + f_x^TV_{xx}f_u + V_x f_{xu} \\ 
    & Q_{ux} = l_{ux} + f_u^TV_{xx}f_x + V_x f_{ux} \\ 
    & Q_{uu} = l_{uu} + f_u^TV_{xx}f_u + V_x f_{uu} \\ 
\end{split}
\end{align}

By taking the derivatives with respect to $\delta u$ and equating to zero, we get a locally linear feedback policy 

\begin{align}\label{eqn:policy}
\begin{split}
    & Q_{uu} \delta u + Q_u + Q_{ux} \delta x = 0 \\
    & \delta u = - Q_{uu}^{-1} [Q_u + Q_{ux} \delta x]
\end{split}
\end{align}

This is equivalent to solving the quadratic approximation of the following optimization problem at time point $k$

\begin{equation}
    \min _{x_k, u_k} l(x_k, u_k) + V_{k + 1}(f(x_k, u_k))
\end{equation}

which results in the following KKT conditions, after dropping the subscript $k$ for simplicity

\begin{equation}    
    \begin{bmatrix}
        l_x + f_x^TV_x\\
        l_u + f_u^TV_x
    \end{bmatrix} = 0
\end{equation}

Because the subproblem is quadratic, we can obtain the solution using one Newton step given by the following set of equations

\begin{equation}\label{eqn:quad_subproblem}
    \begin{bmatrix}
        l_{xx} + f_x^TV_{xx}f_x + V_x f_{xx} & l_{xu} + f_x^TV_{xx}f_u + V_x f_{xu} \\
        l_{ux} + f_u^TV_{xx}f_x + V_x f_{ux} & l_{uu} + f_u^TV_{xx}f_u + V_x f_{uu} 
    \end{bmatrix} \begin{bmatrix}
        \delta x \\
        \delta u
    \end{bmatrix} = - \begin{bmatrix}
        l_x + f_x^TV_x \\
        l_u + f_u^TV_x
    \end{bmatrix}
\end{equation}

However instead of taking a step in both $\delta x$ and $\delta u$ direction, a step is taken only in $\delta u$ direction

\begin{equation}    
    \delta u = - \left [ l_{uu} + f_u^TV_{xx}f_u + V_x f_{uu} \right]^{-1} \left [  l_u + f_u^TV_x + ( l_{ux} + f_u^TV_{xx}f_x + V_x f_{ux} ) \delta x \right]
\end{equation}

Using the definition in equation \ref{eqn:qder}, we arrive at the same equation for $\delta u$ as in equation \ref{eqn:policy}

The equations \ref{eqn:qder} are propagated backward in time starting with the terminal cost and its derivatives. The backward pass gives an update rule for the control inputs as a function of states. The forward pass is then used to get the new state and control sequence for the next iteration. The procedure is repeated until some convergence criteria is reached. 

\begin{align}
\begin{split}
    & V(x_N) = l_N\\
    & V_x(x_N) = l_{N_x}(x_N) \\ 
    & V_{xx}(x_N) = l_{N_{xx}}(x_N) \\ 
\end{split}
\end{align}

\subsection{Constrained Differential Dynamic Programming}

Now, consider a constrained finite time optimal control problem starting at initial state $x_0$

\begin{align}
\label{eqn:opt-constraint}
\begin{split}
    J^*(\textit{X}, \textit{U}) & = \min _{\textit{U}} \sum _{k = 0}^{N - 1} l(x_k, u_k) + l_f(x_N) \\    
    \text{subject to} & \\ 
    & x_{k + 1} = f(x_k, u_k) \quad \forall \ k \in \{ 0, 1, \cdots N - 1 \} \\
    & g_k (x_k, u_k) = 0 \quad \forall \ k \in \{ 0, 1, \cdots N - 1 \} \\
    & h_k (x_k, u_k) \leq 0 \quad \forall \ k \in \{ 0, 1, \cdots N - 1 \} \\
    & g_N (x_N) = 0  \\
    & h_N (x_N) \leq 0 \\
\end{split}
\end{align}

with equality constraints $g_k (x_k, u_k)$ and inequality constraints $h_k (x_k, u_k)$ at every time step $k$.  We add slack variables to convert the inequality constraints to equality constraints and require that the slack variables $s_k$ be positive. Another advantage of adding slack variables is that we can have an infeasible or arbitrary initial trajectory. An interior point algorithm is used to solve a sequence of barrier subproblems that eventually converge to the local optimum of the original unconstrained problem. The barrier subproblem that we solve is as follows 

\begin{align}
\label{eqn:opt-slack}
\begin{split}
    J^*(\textit{X}, \textit{U}, \tau) & = \max _ {\Lambda} \min _{\textit{U}, \textit{M}} \sum _{k = 0}^{N - 1} \left[ l(x_k, u_k) + \lambda _k^T g_k(x_k, u_k) + \frac{\tau}{s_k} (h_k(x_k, u_k) + s_k) - \tau (1^T \log(s_k)) \right ] \\ & + \left[ l_f(x_N) + \lambda _N^T g_N(x_N) + \frac{\tau}{s_N} (h(x_N) + s_N) - \tau (1^T \log(s_N)) \right]\\    
    \text{subject to} & \\ 
    & x_{k + 1} = f(x_k, u_k) \quad \forall \ i \in \{ 0, 1, \cdots N - 1 \} \\
    & s_{k} \geq 0 \quad \forall \ i \in \{ 0, 1, \cdots N \} 
\end{split}
\end{align}

where $\tau$ is the barrier parameter, $\textit{M} = \{s_0, s_1, ..., s_N \}$ are the slack variables and $\Lambda = \{ \lambda _0, \lambda _1, ... \lambda _N \}$ are the lagrange variables corresponding to the equality constraints $g(x_k, u_k)$. We can again apply dynamic programming to this problem by defining the value function as 

\begin{equation}
\label{eqn:val-fun-con}
    V_k(x_k) = \min _{u_k, s_k \geq 0} \max_{\lambda _k} \left[ l(x_k, u_k) + \lambda _k^T g_k(x_k, u_k) + \frac{\tau}{s_k} (h_k(x_k, u_k) + s_k) - \tau (1^T \log(s_k)) + V_{k + 1}(x_{k + 1}) \right]
\end{equation}

To get the DDP update, we solve a similar optimization problem by defining $z_k = [x_k, u_k]$

\begin{equation}
    V_k(x_k) = \min _{z_k, s_k \geq 0} \max_{\lambda _k} \left[ l(x_k, u_k) + \lambda _k^T g_k(x_k, u_k) + \frac{\tau}{s_k} (h_k(x_k, u_k) + s_k) - \tau (1^T \log(s_k)) + V_{k + 1}(x_{k + 1}) \right]
\end{equation}

which has the following KKT conditions, after dropping the subscript k for simplicity

\begin{equation} 
    \begin{bmatrix}
        V_z \\
        V_{\lambda} \\
        V_s
    \end{bmatrix} = 
    \begin{bmatrix}
        l_z + \lambda ^T g_z + \tau S^{-1} h_z + f_z^TV_z\\
        g \\
        h + s 
    \end{bmatrix} = 0
\end{equation}

and requires the following Newton step 

\begin{equation}
    \begin{bmatrix}
        l_{zz} + \lambda ^T g_{zz} + \tau S^{-1} h_{zz} + f_z^TV_{zz}f_z + V_z f & g_z ^T & h_z ^ T S^{-2} \tau\\
        g_z & 0 & 0\\
        h_z & 0 & I
    \end{bmatrix} \begin{bmatrix}
        \delta z \\
        \delta \lambda \\
        \delta s
    \end{bmatrix} = - \begin{bmatrix}
        l_z + \lambda ^T g_z + \tau S^{-1} h_z + f_z^TV_z\\
        g \\
        h + s
    \end{bmatrix}
\end{equation}

\subsection{Regularization}

Let $B_{zz} = l_{zz} + \lambda ^T g_{zz} + \frac{\tau}{s} h_{zz} + f_z^TV_{zz}f_z + V_z f$, $B_{z} = l_z + \lambda ^T g_z + \frac{\tau}{s} h_z + f_z^TV_z $, and $S = diag(s)$. Because we maximize with respect to the lagrange variables $\lambda $, we add a negative definite regularization as a function of the barrier parameter $T = \epsilon (\tau) *I$ to the Hessian matrix. The updated system of equations are

\begin{equation}
    \begin{bmatrix}
        B_{zz} & g_z ^T & - h_z ^ T S^{-2} \tau\\
        g_z & - \mathrm{T} & 0\\
        h_z & 0 & I
    \end{bmatrix} \begin{bmatrix}
        \delta z \\
        \delta \lambda \\
        \delta s
    \end{bmatrix} = - \begin{bmatrix}
        B_z\\
        g \\
        h + s
    \end{bmatrix}
\end{equation}

Eliminating the last two rows gives the following equation for z 

\begin{equation}
    \begin{bmatrix}
        B_{zz} + h_z ^ T \tau S^{-2} h_z + g_z^T \mathrm{T}^{-1} g_z\\
    \end{bmatrix} \begin{bmatrix}
        \delta z \\
    \end{bmatrix} = - \begin{bmatrix} 
        B_z + h_z ^ T \tau S^{-2} h + h_z^T \tau S^{-1} + g_z^T \mathrm{T}^{-1} g\\
    \end{bmatrix}
\end{equation}

and the other variables are given as 

\begin{align}
\begin{split}
    & \delta \lambda = \mathrm{T}^{-1} [ g + g_z \delta z ] \\
    & \delta s = - I^{-1} [h + s + h_z \delta z]
\end{split}
\end{align}

Substituting for $z$ gives

\begin{equation}
    \begin{bmatrix}
        B_{xx} + h_x ^ T \tau S^{-2} h_x + g_x^T \mathrm{T}^{-1} g_x & B_{xu} + h_x^T \tau S^{-2} h_u + g_x^T \mathrm{T}^{-1} g_u\\
        B_{ux} + h_u ^ T \tau S^{-2} h_x + g_u^T \mathrm{T}^{-1} g_x & B_{uu} + h_u^T \tau S^{-2} h_u + g_u^T \mathrm{T}^{-1} g_u\\
    \end{bmatrix} \begin{bmatrix}
        \delta x \\
        \delta u
    \end{bmatrix} = - \begin{bmatrix} 
        B_x + h_x ^ T \tau S^{-2} h + h_x^T \tau S^{-1} + g_x^T \mathrm{T}^{-1} g\\
        B_u + h_u ^ T \tau S^{-2} h + h_u^T \tau S^{-1} + g_u^T \mathrm{T}^{-1} g\\
    \end{bmatrix}
\end{equation}

Comparing the above equation with equation \ref{eqn:quad_subproblem} and equation \ref{eqn:qder}, we get the following values for the quadratic subproblem  

\begin{align}
\label{eqn:qder-con}
\begin{split}
    & \hat{Q} = l + V + \lambda ^T g + \tau S^{-1} (h + s) - \tau 1^T \log(s) \\
    & \hat{Q}_x = l_x + f_x^TV_x + \tau S^{-1}h_x + \lambda ^T g_x + h_x ^ T \tau S^{-2} h + h_x^T \tau S^{-1} + g_x^T \mathrm{T}^{-1} g\\ 
    & \hat{Q}_u = l_u + f_u^TV_x + \tau S^{-1}h_x + \lambda ^T g_x + h_u ^ T \tau S^{-2} h + h_u^T \tau S^{-1} + g_u^T \mathrm{T}^{-1} g\\ 
    & \hat{Q}_{xx} = l_{xx} + f_x^TV_{xx}f_x + V_x f_{xx} + \lambda ^T g_{xx} + \tau S h_{xx} + h_x ^ T \tau S^{-2} h_x + g_x^T \mathrm{T}^{-1} g_x \\ 
    & \hat{Q}_{xu} = l_{xu} + f_x^TV_{xx}f_u + V_x f_{xu} + \lambda ^T g_{xu} + \tau S h_{xu} + h_x ^ T \tau S^{-2} h_u + g_x^T \mathrm{T}^{-1} g_u\\ 
    & \hat{Q}_{ux} = l_{ux} + f_u^TV_{xx}f_x + V_x f_{ux} + \lambda ^T g_{ux} + \tau S h_{ux} + h_u ^ T \tau S^{-2} h_x + g_u^T \mathrm{T}^{-1} g_x\\ 
    & \hat{Q}_{uu} = l_{uu} + f_u^TV_{xx}f_u + V_x f_{uu} + \lambda ^T g_{uu} + \tau S h_{uu} + h_u ^ T \tau S^{-2} h_u + g_u^T \mathrm{T}^{-1} g_u\\ 
\end{split}
\end{align}

To use the control policy equation in \ref{eqn:policy}, $\hat{Q}_{uu}$ should be positive definite. In addition to the global regularization $\mu _1$ proposed in \cite{6386025}, we also add a local regularization $\mu _2$ at every state. The algorithm given in \ref{algo:psd} uses the following equations 

\begin{align}\label{eqn:policy-con}
\begin{split}
    & \tilde{Q}_{xx} = \hat{Q}_{xx} + \mu _2 I \\
    & \tilde{Q}_{ux} = \hat{Q}_{ux} + f_u^T (\mu _1 I ) f_x \\
    & \tilde{Q}_{xu} = \hat{Q}_{xu} + f_x^T (\mu _1 I ) f_u \\
    & \tilde{Q}_{uu} = \hat{Q}_{uu} + f_u^T (\mu _1 I) f_u + \mu _2 I \\  
    & k = - \tilde{Q}^{-1}_{uu}  \hat{Q}_u \\
    & K = - \tilde{Q}^{-1}_{uu}  \hat{Q}_{ux}
\end{split}
\end{align}

\begin{algorithm}
\caption{Backward Pass}
\begin{algorithmic}
    \Require $\textit{X}, \ \textit{U}, \ \textit{M}, \ \Lambda $ trajectory, slack variables and lagrange multipliers respectively, barrier parameter $\tau$ 
    \State $(V_N, V_{N_x}, V_{N_{xx}}) \gets $ Equation \ref{eqn:terminal} \Comment{Terminal cost function approximation}
    \For{$i = N - 1$ to $0$}
    \State $(\hat{Q}_i, \hat{Q}_{x_i}, \hat{Q}_{u_i}, \hat{Q}_{xx_i}, \hat{Q}_{xu_i}, \hat{Q}_{ux_i}, \hat{Q}_{uu_i} ) \gets $ Equation \ref{eqn:qder-con}
    \State $ (\tilde{Q}_{xx_i}, \tilde{Q}_{xu_i}, \tilde{Q}_{ux_i}, \tilde{Q}_{ux_i}) \gets \text{Positive Definite}(\hat{Q}_{xx_i}, \hat{Q}_{xu_i}, \hat{Q}_{ux_i}, \hat{Q}_{ux_i})$
    \State $(k_i, K_i) \gets $ Equation \ref{eqn:policy-con}
    \State $(V_i, V_{x_i}, V_{xx_i}) \gets $ Equation \ref{eqn:state-con}
    \EndFor
    \State $ \mathbb{k} \gets \{ k_0, k_1, \cdots k_{N - 1} \} $
    \State $ \mathbb{K} \gets \{ K_0, K_1, \cdots K_{N - 1} \} $
    \State 
    \Return $\mathbb{k}, \mathbb{K}$
\end{algorithmic}
\label{algo:backward}
\end{algorithm}

\begin{algorithm}
\caption{Positive Definite}
\begin{algorithmic}
    \Require $Q_{xx}, Q_{xu}, Q_{ux}, Q_{uu}$, global regularization parameter $\mu _1$, local regularization parameter $\mu _2$, eigenvalue tolerance $\epsilon _{eig}$
    \State $\mu _2 \gets 0$
    \State $\tilde{Q}_{ux}, \tilde{Q}_{xu}, \tilde{Q}_{uu} \gets $ equation \ref{eqn:policy-con}
    \State $ \tilde{Q}_{zz} \gets \begin{bmatrix}
        \tilde{Q}_{xx} & \tilde{Q}_{xu} \\ \tilde{Q}_{ux} & \tilde{Q}_{uu}
    \end{bmatrix}$ 
    \While{$\min \{ eig(\tilde{Q}_{zz}) \} \geq \epsilon _{eig} $}
    \State $\tilde{Q}_{zz} \gets \tilde{Q}_{zz} + \mu _2 I$
    \State $\mu _2 \gets \mu _2 * 10$
    \EndWhile 
    \State 
    \Return $\tilde{Q}_{xx}, \tilde{Q}_{xu}, \tilde{Q}_{ux}, \tilde{Q}_{uu}$
\end{algorithmic}
\label{algo:psd}
\end{algorithm}

with the same value function updates for intermediate time steps

\begin{align}\label{eqn:state-con}
\begin{split}
    & V = \hat{Q} + \frac{1}{2} k^T \hat{Q}_{uu} k + k^T \hat{Q}_u \\
    & V_{x} = \hat{Q}_{x} + K^T \hat{Q}_{uu} k + K^T \hat{Q}_u + \hat{Q}_{ux}^T k \\
    & V_{xx} = \hat{Q}_{xx} + K^T \hat{Q}_{uu} K + K^T \hat{Q}_{ux} + \hat{Q}_{ux}^T K \\
\end{split}
\end{align}

and value function updates for terminal step

\begin{align}\label{eqn:terminal}
\begin{split}
    V_N & = l_N + \lambda _N^T g_N + \tau S_N^{-1} (h_N + s_N) - \tau 1^T \log(s_N) \\
    V_{N_x} & = l_{x_N} + h_{x_N} ^ T \tau S_N^{-2} h_N + h_{x_N}^T \tau S_N^{-1} + g_{x_N}^T \mathrm{T}^{-1} g_N \\
    V_{N_{xx}} & = l_{xx_N} + \lambda _N^T g_{xx_N} + \tau S_N h_{xx_N} + h_{x_N} ^ T \tau S_N^{-2} h_{x_N} + g_{_N}^T \mathrm{T}^{-1} g_{x_N} \\
\end{split}
\end{align}

\subsection{Line Search}
To ensure convergence, we choose an appropriate step size $\alpha ^j \in (0, 1]$ using the line-search method. We find an upper bound on $\alpha$. Since the slack variables have to be positive at the optimal solution, we ensure that this is maintained at every iterate by solving the following equation proposed by \cite{Wchter2006OnTI} 

\begin{equation}
    \alpha ^{j}_{max} = \max \{ \alpha \in (0, 1] : s + \alpha \delta s \geq (1 - \epsilon) s \}
\end{equation}

where the parameter $\epsilon$ is chosen to be $0.995$. Subsequently, the step size $\alpha ^j \in (0, \alpha ^{j}_{max}]$ for the remaining variables is obtained using backtracking line-search. We accept a step size if there is a sufficient decrease in either the merit function $\phi $ or the constraint violation $\theta$.

\begin{align}
\begin{split}
    & \phi (\mathbb{X}, \mathbb{U}, \mathbb{S}) = \sum _{k = 1}^{N - 1} \left[ l(x_k, u_k) - \tau (1^T \log(s_k)) \right ] + l(x_N) - \tau (1^T \log(s_N)) \\
    & \theta (\mathbb{X}, \mathbb{U}, \mathbb{S}) = \sum _{k = 1}^{N - 1} \left[ || h(x_k, u_k) + s_k ||_1 + || g(x_k, u_k) ||_1 \right] + || h(x_N) + s_N ||_1 + || g(x_N) ||_1 \\
\end{split}
\end{align}

where given a step size $\alpha ^{i} $ and search direction, the remaining variables are updated as follows 

\begin{align}\label{eqn:lagrange-update}
\begin{split}
    u^{j + 1} & = u^{j} + k + \alpha ^{j} K (x^{j + 1} - x^{j})\\
    \lambda ^{j + 1} & = \lambda ^{j} - \alpha ^{j} \mathrm{T}^{-1} g^j - \mathrm{T}^{-1} g^j_{x} (x^{j + 1} - x^{j}) - \mathrm{T}^{-1} g^j_{u} (u^{j + 1} - u^{j}) \\
    s^{j + 1} & = s^{j} - \alpha ^{j} (h^j + s^j) - h^j_{x} (x^{j + 1} - x^{j}) - h^j_{u} (u^{j + 1} - u^{j}) \\
\end{split}
\end{align}

\begin{algorithm}
\caption{Forward Pass}
\begin{algorithmic}
    \Require $\textit{X}^j = \{ x^j_0, x^j_1, \cdots x^j_N \}, \ \textit{U}^j = \{ u^j_0, u^j_1, \cdots u^j_{N - 1} \}, \ \textit{M}^j = \{ s^j_0, s^j_1, \cdots s^j_N \}, \ \Lambda ^j = \{ \lambda ^j_0, \lambda ^j_1, \cdots \lambda ^j_N \}$ trajectory, slack variables and lagrange multipliers at the current iterate, $\alpha ^j$ step size at the current iterate, $ g^j_i, h^j_i $ equality and inequality constraints evaluated at time point $i$ of current iterate $(x_i^j, u_i^j)$, $ g^j_{x_i}, g^j_{u_i}, h^j_{x_i}, h^j_{u_i} $ gradients of equality constraints and inequality constraints with respect to $x$ and $u$ evaluated at time point $i$ of current iterate $(x_i^j, u_i^j)$
    \State $x^{j + 1}_1 \gets x^j_1$ \Comment{Initial condition}
    \For{$i = 0$ to $N - 1$}
    \State $(u^{j + 1}_i, \lambda ^{j + 1}_i, s^{j + 1}_i ) \gets $ equation \ref{eqn:lagrange-update}\
    \State $x^{j + 1}_{i + 1} \gets f(x^{j + 1}_{i}, u^{j + 1}_i)$\
    \EndFor
    \State $(\lambda ^{j + 1}_N, s^{j + 1}_N) \gets $ equation \ref{eqn:lagrange-update} \
    \State 
    \Return $\textit{X}^{j + 1}, \ \textit{U}^{j + 1}, \ \textit{M}^{j + 1}, \ \Lambda ^{j + 1}$ \Comment{Trajectory at next iterate}
\end{algorithmic}
\label{algo:forward}
\end{algorithm} 

\begin{algorithm}
\caption{Backtracking Line-search}
\begin{algorithmic}
    \Require $\textit{X}^j = \{ x^j_0, x^j_1, \cdots x^j_N \}, \ \textit{U}^j = \{ u^j_0, u^j_1, \cdots u^j_{N - 1} \}, \ \textit{M}^j = \{ s^j_0, s^j_1, \cdots s^j_N \}, \ \Lambda ^j = \{ \lambda ^j_0, \lambda ^j_1, \cdots \lambda ^j_N \}, $ trajectory, slack variables and lagrange multipliers at the current iterate, $\alpha ^j$ step size at the current iterate, $\alpha ^j_{max}$ maximum step size, $ \epsilon _{tol} $ tolerance, $\epsilon _{\theta} \in (0, 1)$
    \State $\alpha ^{j} \gets \alpha ^j_{max}$
    \State $(\phi ^j, \theta ^j) \gets (\phi (\textit{X}^j, \textit{U}^j, \textit{M}^j), \ \theta (\textit{X}^j, \textit{U}^j, \textit{M}^j))$
    \While{True}
    \State $(\textit{X}^{j + 1}, \ \textit{U}^{j + 1}, \ \textit{M}^{j + 1}, \ \Lambda ^{j + 1} ) \gets \text{Forward Pass algorithm \ref{algo:forward}} $
    \State $(\phi ^{j + 1}, \theta ^{j + 1} ) \gets (\phi (\textit{X}^{j + 1}, \textit{U}^{j + 1}, \textit{M}^{j + 1}), \ \theta (\textit{X}^{j + 1}, \textit{U}^{j + 1}, \textit{M}^{j + 1}) )$
    \If{$ \theta ^{j + 1} \leq \epsilon _{tol}$}
    \If{$\phi ^{j + 1} \leq \phi ^{j} - \epsilon _{\theta} \theta ^{j}$}
    \State \textbf{Break}
    \EndIf
    \ElsIf{$(\theta ^{j + 1} \leq (1 - \epsilon _{\theta}) \theta ^j ) \ \textbf{or} \ ( \phi ^{j + 1} \leq \phi ^{j} - \epsilon _{\theta} \theta ^{j} )$}
    \State \textbf{Break}
    \EndIf
    \State $\alpha ^{j} \gets 0.5 \alpha ^{j}$
    \EndWhile
    \State
    \Return $\textit{X}^{j + 1}, \ \textit{U}^{j + 1}, \ \textit{M}^{j + 1}, \ \Lambda ^{j + 1}, \alpha ^{j}$
\end{algorithmic}
\label{algo:backtracking}
\end{algorithm}

\subsection{Convergence criteria}
We use the same convergence criteria presented in \cite{Nocedal1997OnTL} for accepting the approximate solution to the barrier problem defined by the barrier parameter $\tau$. 

\begin{equation}\label{eqn:convergence}
     \max \left\{ \theta ^{j + 1}, \ || \hat{Q}_u ||_{\infty} \right\} \leq \epsilon _{\tau} \tau
\end{equation}

And the barrier parameter is updated as follows 

\begin{equation}\label{eqn:tau-update}
    \tau _{j + 1} = \max \left\{ \frac{\epsilon _{tol}}{10}, \min \{ 0.2 \tau _j, \ \tau _j ^{1.5}\} \right\}
\end{equation}

\subsection{Algorithm}

Algorithm \ref{algo:ipddp} gives a pseudo code of the proposed algorithm (https://github.com/siddharth-prabhu/ConstraintDDP). We aim to solve a constrained optimal control problem \ref{eqn:opt-constraint}. We address this by solving a series of barrier subproblems with a given value of the barrier parameter, which is decreased after the convergence criteria for each subproblem are met. Each barrier subproblem is solved using a DDP-like approach. The backward pass of DDP provides the update rule while the forward pass provides the next iterates. 

\begin{algorithm}
\caption{Interior point differential dynamic programming}
\begin{algorithmic}
    \Require Initial control sequence $\mathrm{U}^0$, tolerance $\epsilon _{tol}$
    \State $\textit{X}^0 \gets $ forward simulate dynamics using initial control sequence
    \State $ \textit{M} \gets \max \{h, 1 \}, \ \Lambda \gets 1 ,\ \tau \gets \max \{ \phi (\textit{X}^0, \textit{U}^0, \textit{M}) / N, 10\}$ \Comment{Initialize variables}
    \While{$\tau \geq \epsilon _{tol}$}
    \State search direction $\gets$ Backward Pass algorithm \ref{algo:backward} \Comment{Backward pass}
    \State next iterate $\gets$ Backtracking Line-search algorithm \ref{algo:backtracking} \Comment{Forward pass}
    \If{Equation \ref{eqn:convergence} is True} \Comment{Check convergence}
    \State update $\tau$ using Equation \ref{eqn:tau-update} \Comment{Update barrier parameter}
    \EndIf
    \EndWhile
\end{algorithmic}
\label{algo:ipddp}
\end{algorithm}








\section{Experiments}

\subsection{Inverted Pendulum}
Here we consider the problem of stabilizing an inverted pendulum. The dynamics of the pendulum is as follows

\begin{equation}
    f(x, u) = \begin{bmatrix}
        x_0 + hx_1 \\
        h \sin{x_0} + hu
    \end{bmatrix}
\end{equation}
where the states $ x = [x_0, x_1] $ are the angle and the angular velocity, $u$ is the control input. Step size $h = 0.05$ and initial condition $x = [-\pi, 0]$ is chosen. The control inputs are bounded using inequality constraints $-0.25 \leq u \leq 0.25$ while the terminal states are bounded using equality constraints $x(N) = [0, 0]$. The initial guess for control inputs is randomly chosen from interval $[-0.01, 0.01]$, and the control horizon is chosen to be $N = 500$. The running cost is as follows 

\begin{equation*}
    l(x, u) = 0.025(x_0^2 + x_1^2 + u^2)
\end{equation*}

\begin{figure}[h!]
    \centering
    \begin{subfigure}{.3\textwidth}
      \includegraphics[width = \linewidth, height = 0.15\textheight]{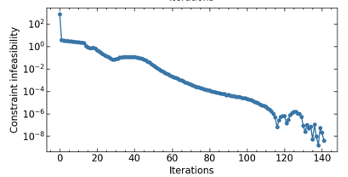}
      \caption{}
      \label{fig:constraints_pendulum}
    \end{subfigure}%
    \centering
    \begin{subfigure}{.3\textwidth}
      \includegraphics[width = \linewidth, height = 0.15\textheight]{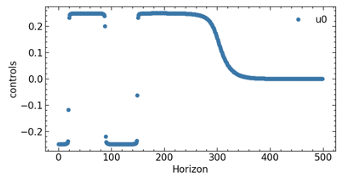}
      \caption{}
      \label{fig:inputs_pendulum}
    \end{subfigure}
    \centering
    \begin{subfigure}{.3\textwidth}
      \includegraphics[width = \linewidth, height = 0.15\textheight]{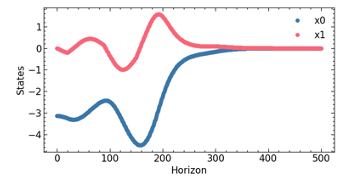}
      \caption{}
      \label{fig:states_pendulum}
    \end{subfigure}
    
    \caption{Solution of the inverted pendulum problem given by algorithm \ref{algo:ipddp}}
    \label{fig:pendulum}
\end{figure}

The algorithm successfully generates the trajectory while obeying the constraints. We plot the maximum constraint violation as the optimization progresses, the optimal control inputs and the optimal states over the control horizon in figure \ref{fig:constraints_pendulum}, \ref{fig:inputs_pendulum}, \ref{fig:states_pendulum} respectively.

\subsection{Continuous Stirred Tank Reactor}
Here, we consider an example of a reactor in which a highly exothermic reaction occurs. A cooling agent maintains the reactor's temperature at a certain value, otherwise, the temperature rapidly increases and degrades the product and/or the reactor. The following equations give the dynamics of the reactor 

\begin{align} \label{eqn:cstr}
    \begin{split}
         \frac{dC}{dt} & = \frac{Q}{V}(C_f - C) - K(T)C\\
         \frac{dT}{dt} &= \frac{Q}{V}(T_f - T) + \frac{-\Delta H}{\rho c_p}K(T)C + \frac{UA}{\rho c_p V}(T_c - T)\\
         \frac{dT_c}{dt} &= \frac{Q_c}{V_c}(T_{cf} - T_c) + \frac{UA}{\rho c_p V_c}(T - T_c)\\
        K(T) & = k_0 e^{\frac{-E_a}{RT}}
    \end{split}
\end{align}

where the states $[C, T, T_c]$ are the concentration of the product in the reactor, the temperature in the reactor, and the temperature of the coolant used respectively. All other constants are as follows - the activation energy $E_a = 72750$ J/gmol, Arrhenius constant $k_0 = 7.2 \times 10^{10}$ 1/min, gas constant $R = 8.314$ J/gmol/K, reactor volume $V = 100$ liters, density $\rho = 1000$ g/liter, heat capacity $c_p = 0.239$ J/g/K, enthalpy of reaction $\Delta H = -5\times 10^4$ J/gmol, heat transfer coefficient $UA = 5 \times 10^4$ J/min/K, feed flowrate $Q = 100$ liters/min, inlet feed concentration $C_f$ gmol/liter, inlet feed temperature $T_f = 300$ K, coolant feed temperature $T_{cf} = 300 $ K, coolant jacket volume $V_c = 20$ liters. The coolant flowrate is the control input $Q_c$. The initial states are given as $C = 0.5, T = 350, T_c = 300$ and the initial steady state value of coolant flowrate is $Q_c = 150$.

\begin{figure}[h!]
    \centering
    \begin{subfigure}{.3\textwidth}
      \includegraphics[width = \linewidth, height = 0.15\textheight]{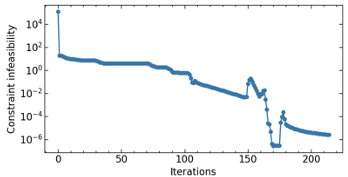}
      \caption{}
      \label{fig:constraints_cstr}
    \end{subfigure}%
    \centering
    \begin{subfigure}{.3\textwidth}
      \includegraphics[width = \linewidth, height = 0.15\textheight]{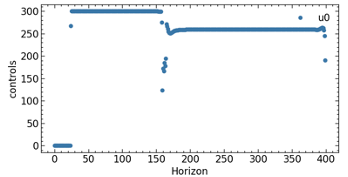}
      \caption{}
      \label{fig:inputs_cstr}
    \end{subfigure}
    \centering
    \begin{subfigure}{.3\textwidth}
      \includegraphics[width = \linewidth, height = 0.15\textheight]{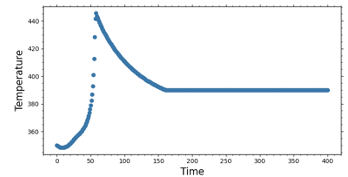}
      \caption{}
      \label{fig:states_cstr}
    \end{subfigure}
    
    \caption{Solution of the continuous stirred tank reactor problem given by algorithm \ref{algo:ipddp}}
    \label{fig:cstr}
\end{figure}

We convert the continuous-time dynamics in \ref{eqn:cstr} into discrete-time using a fixed step $(h = 0.01)$ rk45 solver. We choose a control horizon $N = 400$. The coolant flowrate is used to control the temperature inside the reactor to $390 K$. We bound the control inputs using inequality constraints $\leq u \leq $ and the terminal reactor temperature using equality constraints $T(t = N) = 390$.  We use the following running cost 

\begin{equation*}
    l(x, u) = 50(T - 390)^2
\end{equation*}

We do not penalize the control input in this problem because we don't expect it to go to zero but to rather reach a different steady state. This is observed in the optimal control input and the reactor temperature trajectory resulting from the control input, as shown in figure \ref{fig:inputs_cstr} and \ref{fig:states_cstr}.

\subsection{Car Parking}
\begin{figure}[h!]
    \centering
    \begin{subfigure}{.3\textwidth}
      \includegraphics[width = \linewidth, height = 0.15\textheight]{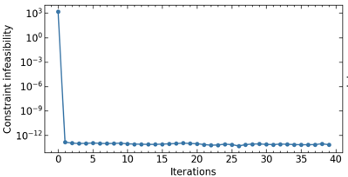}
      \caption{}
      \label{fig:constraints_parking}
    \end{subfigure}%
    \centering
    \begin{subfigure}{.3\textwidth}
      \includegraphics[width = \linewidth, height = 0.15\textheight]{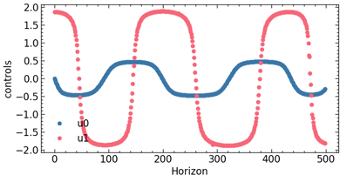}
      \caption{}
      \label{fig:inputs_parking}
    \end{subfigure}
    \centering
    \begin{subfigure}{.3\textwidth}
      \includegraphics[width = \linewidth, height = 0.15\textheight]{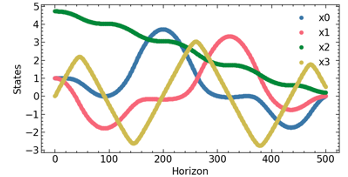}
      \caption{}
      \label{fig:states_parking}
    \end{subfigure}
    
    \caption{Solution of the car parking problem given by algorithm \ref{algo:ipddp}}
    \label{fig:parking}
\end{figure}

We also consider the car parking example from \cite{6907001}. The following equation gives the dynamics of the system

\begin{equation}
    f(x, u) = \begin{bmatrix}
        x_0 + b(x_3, u_0) \cos{x_2} \\
        x_1 + b(x_3, u_0) \sin{x_2} \\
        x_2 + \sin^{-1}\left({\frac{hx_3\sin{u_0}}{d}} \right)\\
        x_3 + h u_1
    \end{bmatrix}
\end{equation}

where the states $x = [x_0, x_1, x_2, x_3]$ are the x-coordinate, the y-coordinates, the cars heading and the velocity of the car, $u = [u_0, u_1]$ are the front wheel steering angle and acceleration. $b(v, w) = d + hv\cos{w} - \sqrt{d^2 - h^2v^2\sin^2{w}}$, $h = 0.03$, and the distance between the front and back axles of a car $d = 2$. The initial condition $x = [1, 1, 3\pi/2, 0]$ is chosen with normally distributed control inputs and control horizon of $N = 500$. The control inputs are bounded with inequality constraints $ -0.5 \leq u_0 \leq 0.5$, $ -2 \leq u_1 \leq 2 $ and the running cost and the terminal cost are as follows 

\begin{align*}
\begin{split}
    & l(x, u) = 0.001H(x_0, 0.1) + 0.001H(x_1, 0.1) + 0.01u_0^2 + 0.0001u_1^2\\ 
    & l_f(x_N) = 0.1H(x_0, 0.01) + 0.1H(x_1, 0.01) + H(x_2, 0.01) + 0.3H(x_3, 1)
\end{split}
\end{align*}

where the function $H(y, z) = \sqrt{y^2 + z^2} - z$. The optimal solution is reached within 40 iterations of the algorithm with a solution accuracy of $10^{-4}$ as shown in figure \ref{fig:parking}

\subsection{Car Obstacle}
Next, we consider a 2D car with the following dynamics 

\begin{equation}
    f(x, u) = \begin{bmatrix}
        x_0 + h x_3 \sin{x_2} \\
        x_1 + hx_3 \cos{x_2} \\
        x_2 + h u_1 x_3 \\
        x_3 + h u_0
    \end{bmatrix}
\end{equation}

where the states $x = [x_0, x_1, x_2, x_3]$ and the control inputs $u = [u_0, u_1]$. Given the initial conditions $x = [0, 0, 0, 0]$, the goal is to reach the terminal state $x(T) = [3, 3, \pi/2 0] $ using control inputs that are bounded using inequality constraints $ \pi / 2 \leq u_0 \leq \pi /2 $, $ -10 \leq u_1 \leq 10 $ and avoiding three obstacles defined using the following inequality constraints 

\begin{align*}
\begin{split}
    0.5^2 - (x_0 - 1)^2 - (x_1 - 1)^2 \leq 0 \\
    0.5^2 - (x_0 - 1)^2 - (x_1 - 2.5)^2 \leq 0 \\
    0.5^2 - (x_0 - 2.5)^2 - (x_1 - 2.5)^2 \leq 0
\end{split}    
\end{align*}

\begin{figure}[h!]
    \centering
    \begin{subfigure}{.3\textwidth}
      \includegraphics[width = \linewidth, height = 0.15\textheight]{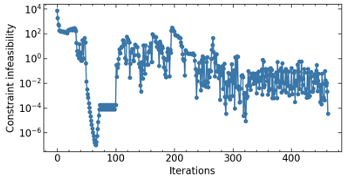}
      \caption{}
      \label{fig:constraints_obstacle}
    \end{subfigure}%
    \centering
    \begin{subfigure}{.3\textwidth}
      \includegraphics[width = \linewidth, height = 0.15\textheight]{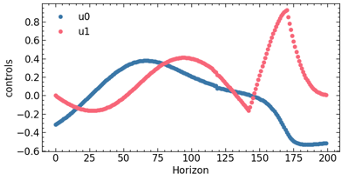}
      \caption{}
      \label{fig:inputs_obstacle}
    \end{subfigure}
    \centering
    \begin{subfigure}{.3\textwidth}
      \includegraphics[width = \linewidth, height = 0.15\textheight]{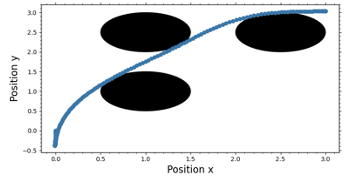}
      \caption{}
      \label{fig:states_obstacle}
    \end{subfigure}
    
    \caption{Solution of the car obstacle problem given by algorithm \ref{algo:ipddp}}
    \label{fig:obstacle}
\end{figure}

The initial guess from control inputs is chosen uniformly between $[-0.01, 0.01]$ and a control horizon of $N = 200$ is chosen. We observe that the algorithm successfully generates a trajectory that avoids the obstacles, as shown in figure \ref{fig:obstacle} 

\section{Conclusion}
Solving DDP with arbitrary equality and inequality state and control constraints has been demonstrated. Adding slack variables to inequality constraints allows us to use an infeasible trajectory as an initial guess. Furthermore, explicit update equations for controls, lagrange variables, and slack variables have been obtained. We also apply this algorithm to a few example systems such as inverted pendulum, continuously stirred tank reactor, car parking and obstacle avoidance. We can easily extend this algorithm in an iLQR setting wherein the second-order derivatives of the dynamics are neglected. A multiple-shooting variant based on \cite{Giftthaler2017AFO, li2023unifiedperspectivemultipleshooting} can also be implemented.

\bibliographystyle{unsrtnat}
\bibliography{references}  

\end{document}